\documentclass[11pt]{article}

\usepackage{amssymb,amsmath,mathrsfs,mathtools,amsthm,latexsym,hyperref}
\usepackage[a4paper,includeheadfoot,left=30mm,right=30mm,top=25mm,bottom=30mm]{geometry} 
\usepackage{verbatim}
\usepackage{color}

\newtheorem{thm}{Theorem}[section]

\theoremstyle{definition}

\theoremstyle{remark}
\newtheorem{rmk}[thm]{Remark}

\newcommand{\hsforall}{\hspace{1mm}\forall\hspace{1mm}}						 
\newcommand{\hsexists}{\hspace{1mm}\exists\hspace{1mm}} 					 
\renewcommand{\Re}{\operatorname*{Re}}                             
\renewcommand{\d}{\ensuremath{\,\mathrm{d}}}							         
\renewcommand{\leq}{\leqslant}                                     
\newcommand{\BE}{\begin{equation}}                                 
\newcommand{\EE}{\end{equation}}                                   
\newcommand{\BES}{\begin{equation*}}                               
\newcommand{\EES}{\end{equation*}}                                 
\newcommand{\BP}{\begin{pmatrix}}                                  
\newcommand{\EP}{\end{pmatrix}}                                    
\newcommand{\N}{\mathbb{N}}                                        
\newcommand{\R}{\mathbb{R}}                                        
\newcommand{\C}{\mathbb{C}}                                        

\usepackage{color, graphicx, graphics}
\newcommand{\ba}{\begin{array}}
\newcommand{\ea}{\end{array}}
\newcommand{\ud}{\,\mathrm{d}}
\newcommand{\bea}{\begin{eqnarray}}
\newcommand{\eea}{\end{eqnarray}}
\graphicspath{{./gfx/}}

\newcommand{\superscript}[1]{\ensuremath{^{\textrm{#1}}}}

\newcommand{\Thns}[0]{\superscript{th}}
\newcommand{\Th}[0]{\Thns~}
\newcommand{\stns}[0]{\superscript{st}}
\newcommand{\st}[0]{\stns~}

\def\clap#1{\hbox to 0pt{\hss#1\hss}}

\numberwithin{equation}{section}

\title{Heat equation on a network using the Fokas method}
\author{N. E. Sheils$^1$ and D. A. Smith$^2$\\
\footnotesize 1. Department of Applied Mathematics, University of Washington, WA \\
\footnotesize 2. Department of Mathematical Sciences, University of Cincinnati, OH \\
\footnotesize Current address: Department of Mathematics, University of Michigan, MI \\
\footnotesize email\textup{: \texttt{daasmith@umich.edu}}
}

\begin{document}
\maketitle

\begin{abstract}
The problem of heat conduction on networks of multiply connected rods is solved by providing an explicit solution of the one-dimensional heat equation in each domain. The size and connectivity of the rods is known, but neither temperature nor heat flux are prescribed at the interface. Instead, the physical assumptions of continuity at the interfaces are the only conditions imposed. This work generalizes that of Deconinck, Pelloni, and Sheils~\cite{DeconinckPelloniSheils} for heat conduction on a series of one-dimensional rods connected end-to-end to the case of general configurations.
\end{abstract}

\section{Introduction}

Suppose $V$ is a set of vertices, and $E$ is a set of directed edges between those vertices so that $(V,E)$ is a finite connected directed graph. Associated with each edge $r\in E$ is the length $L_r\in(0,\infty]$. Suppose additionally if $L_r=\infty$, then the vertex at which $r$ terminates has no other edges. For each edge $r\in E$, we define the open interval $\Omega_r=(0,L_r)$.  We consider the graph to represent a configuration of narrow metal rods joined at the vertices in perfect thermal contact. In this paper we study the distribution of heat on such a configuration. 

In order to derive the proper equations and interface conditions, we generalize what is outlined in~\cite{Kevorkian} for problems with multiple rods.  Assume each edge $r\in E$ is a rod made of some heat-conducting material with density $\rho_r(x)$ and unit cross-sectional area.  We assume the surface of the rod is perfectly insulated so no heat is lost or gained through this surface.  If we consider an infinitesimal section of length $\mathrm{d} x$ for $x\in\Omega_r$, then $\mathrm{d} C^{r}$, the heat content in the section, is proportional to the mass and the temperature, $q^{r}(x,t)$. That is
\BES
	\ud C^{r}(t)=c_r(x)\rho_r(x)q^{r}(x,t) \ud x,
\EES
where $c_r(x)$ is the specific heat on $r$ and $x\in\Omega_r$. Thus, the total heat content in the interval $x_1\leq x\leq x_2$ where $x_1,x_2\in\Omega_r$ is
\BES
	C^{r}(t)=\int_{x_1}^{x_2}c_r(x)\rho_r(x)q^{r}(x,t)\ud x.
\EES

Fourier's Law for heat conduction~\cite{Fourier} states that the rate of heat flowing into a body is proportional to the area of that element and to the outward normal derivative of the temperature at that location. The constant of proportionality is the thermal conductivity, $k_r(x)$.  In our example, the net inflow of heat through the boundaries $x_1$ and $x_2$ is
\BE \label{defn_r}
	R(t)=k_r(x_2)q^{r}_x(x_2,t)-k_r(x_1)q^{r}_x(x_1,t),
\EE
where $x_1,x_2\in\Omega_r$. Conservation of heat implies $\frac{\ud}{\ud t}C^{r}(t)=R(t)$ along each rod.  That is,
\BE\label{conserve_r}
	\frac{\ud}{\ud t} \int_{x_1}^{x_2}c_r(x) \rho_r(x)q^{r}(x,t)\ud x=k_r(x_2)q^{r}_x(x_2,t)-k_r(x_1)q^{r}_x(x_1,t),
\EE
for every $r$. This is a typical conservation law. On each rod we assume constant material properties.  That is $c_r(x)=c_r, \rho_r(x)=\rho_r$, and $k_r(x)= k_r$ for $x\in\Omega_r$ and $0<t<T$. Expressing $R(t)$ as the integral of a derivative and rewriting~\eqref{conserve_r}, we have
\BE
	\int_{x_1}^{x_2}\left(c_r \rho_r q^{r}_t(x,t)\ud x-\frac{\partial}{\partial x}\left(k_r q^{r}_x(x,t)\right)\right) \ud x=0,
\EE
for $x_1,x_2\in\Omega_r$.  Since this is true for any $x_1$ and $x_2$ in this range, it follows that the integrand must vanish.  That is,
$$q^{r}_t(x,t)-\frac{k_r}{c_r\rho_r}q^{r}_{xx}(x,t)=0.$$

Next, we scale $x$ such that $\hat{x}=\frac{x}{c_r\rho_r}$ for $x\in\Omega_r$.  Note that such a scaling also affects $L_r$ and $\Omega_r$ and in what follows we assume all quantities are properly scaled.  Dropping the~$\hat{\cdot}$, and defining $\sigma_r^2=k_r c_r\rho_r$ as the (scaled) thermal diffusivity in rod $r$ we have 
\begin{subequations}
\begin{align} \label{eqn:PDE}
q^{r}_t&=\sigma_r^2 q^{r}_{xx},  &x&\in\Omega_r,~t>0, &\quad r\in E, \\ \label{eqn:InitCond}
	q^{r}(x,0) &= q_0^{r}(x), & x &\in \overline{\Omega_r}, &\quad r\in E,
\end{align}
\end{subequations}
where $q_0^{r}\in\mathcal{S}\left(\overline{\Omega_r}\right)$ is  given initial data for each $r\in E$ with $\mathcal{S}\left(\overline{\Omega}\right)$ a Schwartz space of smooth, rapidly decaying functions restricted to the closure of the interval $\Omega$.  Wherever $\Omega$ is finite, $\mathcal{S}\left(\overline{\Omega}\right)=C^\infty\left(\overline{\Omega}\right)$.  

We also assume that the rods at each vertex $v$ are in perfect thermal contact~\cite{CarslawJaeger}.  The temperature at $x=0$ of the rods emanating from $v$, $q^{r}(0,t)$, is the same as the temperature at $x=L_r$ of the rods terminating at $v$, $q^r(L_r,t)$.   That is, 
\begin{equation}\label{cont_q}
\begin{split}
\begin{aligned}
		\hsexists \theta_v:[0,T]\to\R : &\hsforall r\in E, \mbox{ if } r \mbox{ emanates from } v \mbox{ then } q^{r}(0,t) = \theta_v(t), & t &\in[0,T], \\
		&\hsforall r\in E, \mbox{ if } r \mbox{ terminates at } v \mbox{ then } q^{r}(L_r,t) = \theta_v(t), & t &\in[0,T].
\end{aligned}
\end{split}
\end{equation}
Note that the functions $\theta_v$ are not data of the problem; they appear in the interface condition only for notational convenience. If vertex $v$ has $p$ edges then~\eqref{cont_q} gives $p-1$ conditions.  

In order to find a condition on the set of spatial derivatives of $q^r$ we require an appropriately scaled~\eqref{conserve_r} valid on a region centered at the vertex $v$. Consider a ball centered at a vertex $v$, with radius $\epsilon>0$ sufficiently small so that all rods are straight lines within the ball and no other rods intersect the ball. Then, similar to Equation~\eqref{defn_r} we have
\BE \label{defn2_r}
	R(t)=\sum_{\substack{s\in E: \\ s\text{ emanates} \\ \text{from } v}} \sigma_s^2(\epsilon)q^{s}_x(\epsilon,t) - \sum_{\substack{r\in E: \\ r\text{ terminates} \\ \text{at } v}} \sigma_r^2(L_r-\epsilon)q^{r}_x(L_r-\epsilon,t).
\EE
By conservation of energy,
\begin{equation}\label{conserve2}
\begin{split}
	\frac{\ud}{\ud t}\left( \sum_{\substack{s\in E:\\ \text{$s$ emanates}\\\text{from $v$}}} c_s^2\rho_s^2 \lim_{\epsilon\to0}\int_{0}^{\epsilon}q^s(x,t)\ud x+ \sum_{\substack{r\in E:\\ \text{$r$ terminates}\\\text{at $v$}}} c_r^2\rho_r^2\lim_{\epsilon\to0}\int_{L_r-\epsilon}^{L_r}q^r(x,t)\ud x\right) \\
	=\lim_{\epsilon\to0}\sum_{\substack{s\in E:\\ \text{$s$ emanates}\\\text{from $v$}}} \sigma_s^2 q^s_x(\epsilon,t)-\lim_{\epsilon\to0}\sum_{\substack{r\in E:\\ \text{$r$ terminates}\\\text{at $v$}}} \sigma_r^2q^r_x(L_r-\epsilon,t).
\end{split}
\end{equation}
The left-hand-side of~\eqref{conserve2} is zero by~\eqref{cont_q}. This implies
\begin{equation}\label{cont_qx}
	\sum_{\substack{s\in E:\\ \text{$s$ emanates}\\\text{from $v$}}} \sigma_s^2q^s_x(0,t)-\sum_{\substack{r\in E:\\ \text{$r$ terminates}\\\text{at $v$}}} \sigma_r^2q^r_x(L_r,t)=0.
\end{equation}
Thus the appropriate sum of the heat flux is continuous across the interface.

For each finite endpoint vertex $v\in V$ which has only one finite-length edge $r_v\in E$ (and no infinite-length edges) we prescribe a general Robin boundary condition 
\begin{align}\label{RobinBC}
	\beta_0^{v}q^{r_v}(X,t)+\beta_1^{v}\partial_x q^{r_v}(X,t) &= f_{v}(t), & t &\in[0,T],
\end{align}
where $f_{r_v}\in C^\infty[0,T]$ is known boundary datum and $\beta_0^v,\beta_1^v\in\R$ where $X=0$ if $r_v$ emanates from $v$, and $X=L_{r_v}$ if $r_v$ terminates at $v$. If the graph has $m_f$ finite edges and $m_i$ infinite edges, then~\eqref{cont_q},~\eqref{cont_qx}, and~\eqref{RobinBC} together prescribe a total of $2m_f+m_i$ boundary and interface conditions. 

We insist that the initial data are mutually compatible according to the interface conditions. That is, we require for each interface vertex $v$,
\begin{gather*}
	\sum_{\substack{s\in E:\\s \text{ emanates}\\\text{from }v}}\sigma_s^2\partial_xq^{s}_0(0) - \sum_{\substack{r\in E:\\r \text{ terminates}\\\text{at }v}}\sigma_r^2\partial_xq^{r}_0(L_r) = 0, \\
	\hsforall s\in E, \mbox{ if } s \mbox{ emanates from } v \mbox{ then } q^{s}_0(0) = \theta_v(0), \\
	\hsforall r\in E, \mbox{ if } r \mbox{ terminates at } v \mbox{ then } q^{r}_0(L_r) = \theta_v(0).
\end{gather*}
We also require the initial and boundary conditions are compatible in the sense that for each finite endpoint vertex $v$ with edge $r_v$,
\BES
	\beta_0^vq^{r_v}_0(X)+\beta_1^{v}\partial_x q^{r_v}_0(X) = f_{v}(0).
\EES

Deconinck, Pelloni, and Sheils~\cite{DeconinckPelloniSheils} apply the Fokas method, alternatively called the unified transform method, to a subclass of the problems described above in which each vertex has at most two edges: a configuration of metal rods with different diffusivities joined end-to-end along a line. They derive the so-called \emph{global relation} for each rod using Green's theorem applied to the space-time domain, and follow the approach of~\cite{DTV2014a} to obtain an implicit integral representation for the solution, known as the Ehrenpreis form~\cite{Fok2000a}.

This representation depends upon (time-transforms of) all boundary and interface values $\partial_x^jq^r(0,t)$, $\partial_x^jq^r(L_r,t)$ for $j\in\{0,1\}$, $r\in E$, $t\in[0,T]$, so cannot be considered an effective solution representation. Using the boundary and interface conditions together with the global relation, the authors are able to formulate and solve a system of linear equations for these transformed boundary values in a number of particular cases: two or three semi-infinite, or finite rods. It is clear that the method used would work for a configuration with more than three rods, but the linear algebra becomes cumbersome to solve by hand.

The principal contribution of the present work is to show that the Fokas method may also be applied to the more general configurations of rods described above. The heat equation on graphs is important in engineering applications such as building and analyzing heat exchanger networks~\cite{ShivakumarNarasimhan, YeeGrossmann} and ``channelling'' heat through materials of differing conductivity~\cite{Bejan}.

\subsection{Organization of paper}

In Section~\ref{sec:Implicit} we derive a global relation on each domain $\Omega_r\times(0,T)$ and use it to obtain in implicit integral representation for the solution in Ehrenpreis form.

In the proceeding sections we present explicit solutions for several cases of relatively simple configurations of rods.  Specifically, in Section~\ref{sec:mSemiInf} we study the problem of arbitrarily many semi-infinite rods joined at a single vertex. In Section~\ref{sec:mParallel} we consider arbitrarily many finite parallel rods between a pair of common vertices, and in Section~\ref{sec:mFin} we analyze arbitrarily many finite rods emanating from a single common vertex, and present an example of such a network of three rods.

We conclude with some remarks on the applicability of the method presented herein to related problems.

\section{Implicit integral representation of the solution} \label{sec:Implicit}

\subsection{Global relation} \label{sec:Implicit:GR}

The so-called \emph{global relation} is an essential tool in the Fokas method. As our problem is posed on domains $\Omega_r$ for each $r\in E$ we will derive $|E|$ global relations. Each global relation links the boundary and interface values with the initial value for a particular rod, and holds for all values of the complex spectral parameter $\lambda$ within a half or full plane. In contrast with~\cite{DeconinckPelloniSheils}, but following the innovation of~\cite{Asvestas, Mantzavinos}, we scale $\lambda$ in a way that will simplify later calculations.

For $(x,t)\in\Omega_r$ and $\lambda\in\C$, we define the functions
\begin{align*}
	X_r(x,t,\lambda) &= e^{-ix\lambda/\sigma_r + \lambda^2t} q^r(x,t), \\
	Y_r(x,t,\lambda) &= e^{-ix\lambda/\sigma_r + \lambda^2t} \left(\sigma_ri\lambda + \sigma_r^2\partial_x\right)q^r(x,t),
\end{align*}
which have the property $(\partial_tX_r-\partial_xY_r)(x,t,\lambda)=0$ if and only if  $q^r$ satisfies equation~\eqref{eqn:PDE}.  Suppose that $q^r$ satisfies the partial differential equation and the initial condition~\eqref{eqn:InitCond}. Then Green's theorem~\cite{AblowitzFokas} applied to $\Omega_r(t)=(0,L_r)\times(0,t)$ yields
\BE \label{eqn:GR.Derivation.1}
	0 = \int_{\Omega_r(t)} (\partial_tX-\partial_xY)(x,s,\lambda) \ud x \ud s = \int_{\partial \Omega_r(t)} (Y\ud t + X\ud x).
\EE
Rearranging equation~\eqref{eqn:GR.Derivation.1} and using the notation
\begin{align}
	\hat{q}_0^r(\lambda) &= \int_0^{L_r} e^{-ix\lambda} q_0(x) \ud x, & \lambda &\in \C, \\
	\hat{q}^r(\lambda;t) &= \int_0^{L_r} e^{-ix\lambda} q(x,t) \ud x, & \lambda &\in \C, & t &\in[0,T], \\
	g^{r}_0(\omega, t) &= \int_0^t e^{\omega s}q^{r}(0,s) \ud s, & \omega &\in \C, & t &\in[0,T], \\
	g^{r}_1(\omega, t) &= \int_0^t e^{\omega s}q^{r}_x(0,s) \ud s, & \omega &\in \C, & t &\in[0,T], \\
	h^{r}_0(\omega, t) &= \int_0^t e^{\omega s}q^{r}(L_r,s) \ud s, & \omega &\in \C, & t &\in[0,T], \\
	h^{r}_1(\omega, t) &= \int_0^t e^{\omega s}q^{r}_x(L_r,s) \ud s, & \omega &\in \C, & t &\in[0,T],
\end{align}
yields the global relation
\begin{subequations} \label{eqn:GR.t}
\begin{equation} \label{eqn:GR.Fin.t}
\begin{split}
\hat{q}_0^{r}\left(\frac{\lambda}{\sigma_r}\right) - e^{\lambda^2t}\hat{q}^{r}\left(\frac{\lambda}{\sigma_r};t\right)={} &\sigma_r i\lambda g_0^{r}(\lambda^2,t) + \sigma_r^2 g_1^{r}(\lambda^2,t) \\
&- e^{-i\lambda L_r/\sigma_r} \left( \sigma_r i\lambda h_0^{r}(\lambda^2,t) + \sigma_r^2 h_1^{r}(\lambda^2,t) \right),
\end{split}
\end{equation}
for $\lambda\in\C$ and $t\in[0,T]$ which holds for each $r\in E$ such that $L_r<\infty$. Similarly,
\BE \label{eqn:GR.Inf.t}
	\sigma_r i\lambda g_0^{r}(\lambda^2,t) + \sigma_r^2 g_1^{r}(\lambda^2,t) = \hat{q}_0^{r}\left(\frac{\lambda}{\sigma_r}\right) - e^{\lambda^2t}\hat{q}^{r}\left(\frac{\lambda}{\sigma_r};t\right), \quad \lambda\in\C^-, \; t\in[0,T],
\EE
\end{subequations}
for each $r\in E$ such that $L_r=\infty$.

\subsection{Implicit integral representation of the solution}

\begin{figure}[htbp]
\begin{center}
	\includegraphics{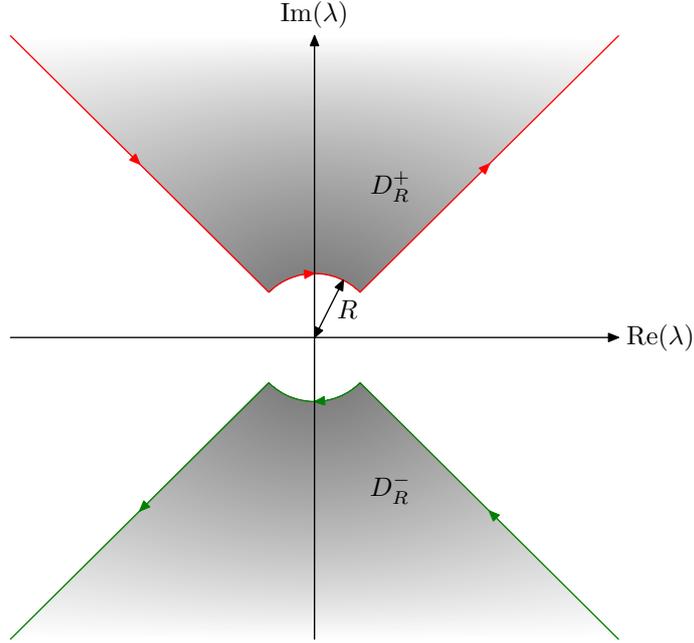}
	\caption{The sectors $D^+_R$ and $D^-_R$, and their oriented boundaries.  \label{fig:heat_DRpm}}
	\end{center}
\end{figure}

Rearranging the global relation~\eqref{eqn:GR.Inf.t}, and taking inverse Fourier transforms in~\eqref{eqn:GR.Fin.t} and~\eqref{eqn:GR.Inf.t} respectively we have
\begin{subequations}
\begin{equation} \label{eqn:ImplicitSol.Derivation.1.Fin}
\begin{split}
q^{r}(x,t) ={} &\frac{1}{2\pi} \int_{-\infty}^\infty e^{i\lambda  x - \sigma_r^2 \lambda^2t} \hat{q}_0^{r}(\lambda) \ud \lambda \\
	&- \frac{1}{2\pi} \int_{-\infty}^\infty e^{ix\lambda / \sigma_r -\lambda ^2t} \left(i\lambda  g_0^{r }(\lambda ^2, t) + \sigma_r g_1^{r}(\lambda ^2, t) \right) \ud \lambda \\
	&+ \frac{1}{2\pi} \int_{-\infty}^\infty e^{i\lambda /\sigma_r(x-L_r) -\lambda ^2t} \left(i\lambda  h_0^{r}(\lambda ^2, t) + \sigma_r h_1^{r}(\lambda ^2, t) \right) \ud \lambda ,
	\end{split}
\end{equation}
and
\begin{equation} \label{eqn:ImplicitSol.Derivation.1.Inf}
\begin{split}
	q^{r}(x,t) ={} &\frac{1}{2\pi} \int_{-\infty}^\infty e^{i\lambda  x - \sigma_r^2\lambda ^2t} \hat{q}_0^{r}(\lambda ) \ud \lambda  \\
	&- \frac{1}{2\pi} \int_{-\infty}^\infty e^{ix\lambda / \sigma_r -\lambda ^2t} \left(i\lambda  g_0^{r}(\lambda ^2, t) + \sigma_r g_1^{r}(\lambda ^2, t) \right) \ud \lambda .
\end{split}
\end{equation}
\end{subequations}
We define the sectors
\begin{figure}[htbp]
	\begin{center}
	\includegraphics{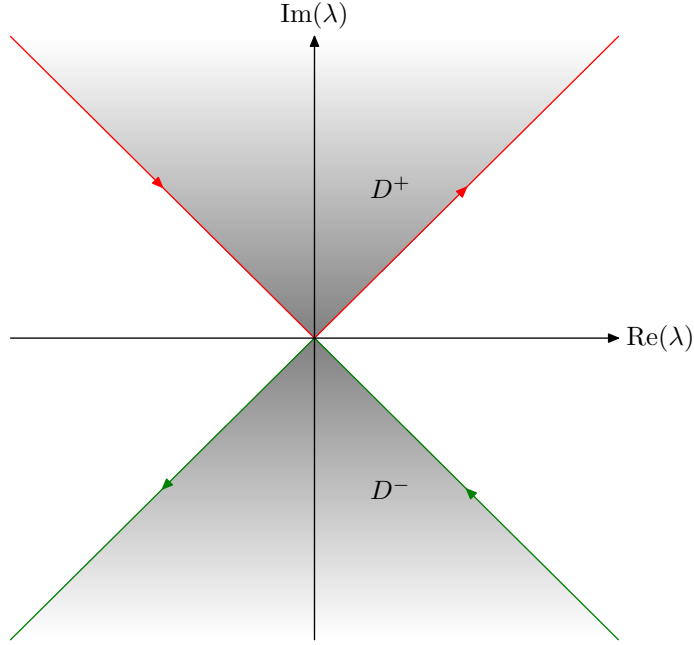}
	\caption{The sectors $D^+$ and $D^-$, and their oriented boundaries.   \label{fig:heat_Dpm}}
	\end{center}
\end{figure}
\BE \label{eqn:Defn.D}
	D^\pm = \{ \lambda\in\C : \tfrac{\pi}{4} < \arg(\pm\lambda) < \tfrac{3\pi}{4} \},
\EE
as shown in Figure~\ref{fig:heat_Dpm}.  The integrand of the second integral in~\eqref{eqn:ImplicitSol.Derivation.1.Inf} is entire and decays as $\lambda\to\infty$ from within $\C^+\setminus D^+$. Hence, by Jordan's Lemma~\cite{AblowitzFokas} we can replace the contour of integration of the second integral by $\int_{\partial D^+}$. Further, by analyticity we make an additional contour deformation from $\partial D^\pm$ to $\partial D^\pm_R$, where
\BE
	D^\pm_R = \{ \lambda\in D^\pm : |\lambda| > R \},
\EE
and $R>0$ is an arbitrary constant. An appropriate (sufficiently large) value of this constant may be chosen for any individual problem as in Figure~\ref{fig:heat_DRpm}.

We have arrived at the implicit integral representation of the solution $q^{r}$ known as the Ehrenpreis form. For each $r\in E$ such that $L_R=\infty$,
\begin{subequations} \label{eqn:ImplicitSol.t}
\begin{equation} \label{eqn:ImplicitSol.Inf.t}
\begin{split}
	q^{r}(x,t) ={} &\frac{1}{2\pi} \int_{-\infty}^\infty e^{i\lambda x - \sigma_r^2\lambda^2t} \hat{q}_0^{r}(\lambda) \d \lambda \\
	&- \frac{1}{2\pi} \int_{\partial D^+_R} e^{ix\lambda/ \sigma_r -\lambda^2t} \left(i\lambda g_0^{r}(\lambda^2,t) + \sigma_r g_1^{r}(\lambda^2,t) \right) \d \lambda.
\end{split}
\end{equation}
Similarly, for each $r\in E$ such that $L_r<\infty$,
\begin{equation} \label{eqn:ImplicitSol.Fin.t}
\begin{split}
	q^{r}(x,t) ={} &\frac{1}{2\pi}\int_{-\infty}^\infty e^{i\lambda x - \sigma_r^2\lambda^2t} \hat{q}_0^{r}(\lambda) \d \lambda \\
	&- \frac{1}{2\pi} \int_{\partial D_R^+} e^{ix\lambda/\sigma_r - \lambda^2t} \left(i\lambda g_0^{r}(\lambda^2,t) + \sigma_r g_1^{r}(\lambda^2,t) \right) \d \lambda \\
	&- \frac{1}{2\pi} \int_{\partial D_R^-} e^{i (x-L_r)\lambda/\sigma_r - \lambda^2t} \left(i\lambda h_0^{r}(\lambda^2,t) + \sigma_r h_1^{r}(\lambda^2,t) \right) \d \lambda.
\end{split}
\end{equation}
\end{subequations}

The process of using the initial, boundary, and interface conditions to determine the $t$-transformed boundary values is known as a generalized spectral Dirichlet-to-Neumann map; generalized because the map is from the known boundary and initial data to the unknown boundary and interface values, rather than from 0\Th order boundary data to 1\st order boundary values; spectral because the map is formulated in spectral space instead of coordinate space. In each of Sections~\ref{sec:mSemiInf}--\ref{sec:mFin}, we implement the generalized spectral Dirichlet-to-Neumann map for particular classes of configurations of rods.

\section{$m$ semi-infinite rods} \label{sec:mSemiInf}
\begin{figure}
	\begin{center}
	\includegraphics{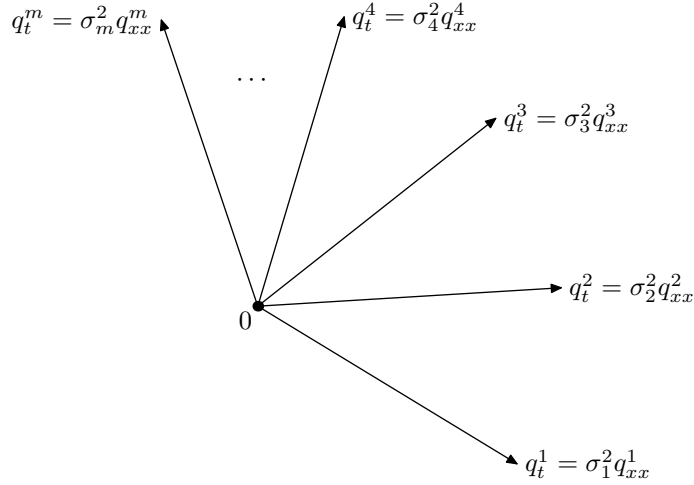}
	\caption{$m$ semi-infinite rods.}
	\label{fig:mSemiInf}
	\end{center}
\end{figure}

In this section, we consider the case of $m\in\N$ semi-infinite rods of differing thermal diffusivity, joined at a single vertex. The graph is as shown in Figure~\ref{fig:mSemiInf}. Of the $m+1$ vertices, $m$ are each the terminus of a single edge, and all $m$ edges emanate from the other vertex. As $|E|=m$, we find it notationally convenient to identify the set of edges $E$ with the set $\{1,2,\ldots,m\}$. The single interface vertex is denoted $v$.

In what follows we find explicit representations of the $t$-transformed interface values appearing in Equation~\eqref{eqn:ImplicitSol.Inf.t} in terms of only the initial data and certain remainder terms. We then show the remainder terms do not to contribute to $q^{r}$.

For each $1\leq r\leq m$, we evaluate the global relation~\eqref{eqn:GR.Inf.t} at $\lambda\mapsto-\lambda$ to obtain the relations
\BE \label{eqn:neg.GR.Inf.t}
	-\sigma_r i\lambda g_0^{r}(\lambda^2,t) + \sigma_r^2 g_1^{r}(\lambda^2,t) = \hat{q}_0^{r}\left(\frac{-\lambda}{\sigma_r}\right) - e^{\lambda^2t}\hat{q}^{r}\left(\frac{-\lambda}{\sigma_r};t\right),
\EE
for all $\lambda\in\C^+$, where we have used that $g_j^{r}(\lambda^2,t)$ is invariant under $\lambda\to-\lambda$. This is essential because we require equations for $g_j^r(\lambda^2,t)$ valid for $\lambda\in D^+\subset\C^+$, whereas the original global relation equation~\eqref{eqn:GR.Inf.t} was valid only for $\lambda\in\C^-$. Equation~\eqref{eqn:neg.GR.Inf.t} represents a system of $m$ linear equations in the $2m$ unknown, $t$-transformed interface values $g_j^{r}(\lambda^2,t)$, for $1\leq r\leq m$, $j\in\{0,1\}$.  

Appropriately transforming the interface conditions~\eqref{cont_q} and~\eqref{cont_qx} for this problem we find
\begin{align}
	&g_0(\lambda^2,t) := g_0^{r}(\lambda^2,t)=\int_0^te^{\lambda^2s}\theta_v(s)\ud s,~~~1\leq r\leq m,\label{msemiinf.IC}\\ 
	&\sum_{r=1}^m \sigma_r^2 g_1^{r}(\lambda^2,t)= 0, \label{eqn:msemiinf.Neumann}
\end{align}
for all $\lambda\in\C$. Equation~\eqref{msemiinf.IC} reduces the number of unknown functions from $2m$ to $m+1$.  Equations~\eqref{eqn:neg.GR.Inf.t},~\eqref{eqn:msemiinf.Neumann} provide $m+1$ linear equations in these unknown functions.  We represent this system as a matrix:
\BE
	\mathcal{A} X=Y,
\EE
where 
\begin{align}
	X &= \left( g_0(\lambda^2,t), \sigma_1^2 g_1^{1}(\lambda^2,t),\ldots,\sigma_m^2g_1^{m}(\lambda^2,t) \right)^\top, \\
	Y &= \left( 0, \hat{q}_0^1\left(\frac{-\lambda}{\sigma_1}\right) - e^{\lambda^2t} \hat{q}^1\left(\frac{-\lambda}{\sigma_1};t\right),\ldots,\hat{q}_0^{m}\left(\frac{-\lambda}{\sigma_m}\right) - e^{\lambda^2t} \hat{q}^{m}\left(\frac{-\lambda}{\sigma_m}; t\right)\right)^\top, \\
	\mathcal{A} &= \BP 0 & 1 & 1 & \cdots & 1 \\ -i\lambda\sigma_1 & 1 & 0 & \cdots & 0 \\ -i\lambda\sigma_2 & 0 & 1 & \cdots & 0 \\ \vdots & \vdots & \vdots & \ddots & \vdots \\ -i\lambda\sigma_m & 0 & 0 & \cdots & 1 \EP.
\end{align}
For notational convenience, the entries of $Y$ are denoted $Y_r$, where $r$ is enumerated $0,\ldots,m$.  Similarly, the rows of the other matrices are counted from $0$ to $m$.

Clearly, the ``data'' vector $Y$ depends not only upon the initial data, but also upon the Fourier transforms of the solutions at time $t$ which are not given data of the problem. Nevertheless, we proceed to solve the system as if $Y$ is composed purely of known data and show afterwards that terms involving $\hat{q}^r(\cdot;t)$ make no contribution to the solution representation.

Subtracting the sum of rows $1,\ldots,m$ from row $0$ in matrix $\mathcal{A}$ it is immediate that
$$
	\det\mathcal{A} = i\lambda\sum_{r=1}^m\sigma_r.
$$
Let $\mathcal{A}_j$ be the matrix formed by replacing column $j$ of $\mathcal{A}$ with the column vector $Y$. By Cramer's rule,
$$
	g_0(\lambda^2,t)=\frac{\det \mathcal{A}_0}{\det \mathcal{A}},
	\qquad
	g_1^{r}(\lambda^2,t)=\frac{\det \mathcal{A}_{r}}{\sigma_r^2 \det \mathcal{A}}, \qquad 1\leq r\leq m.
$$
Trivially,
$$
	\det\mathcal{A}_0 = -\sum_{p=1}^mY_p.
$$
To find $\det\mathcal{A}_r$ we subtract the $j$\Th row for $j=1,2,\ldots,r-1,r+1,\ldots,m$ from the first row. Then, for each $j=r$, $j=r-1$, $\ldots$ , $j=2$, we switch the $j$\Th and $(j-1)$\st rows and switch the $j$\Th and $(j-1)$\st columns. The resulting matrix is
\BE
	\mathcal{B}_r = 
	\BP
		i\lambda \displaystyle\sum_{\substack{p=1\\ p \neq r \\\text{~}}}^m\sigma_p & -\displaystyle\sum_{\substack{p=1\\ p \neq r}}^m Y_p & 0 & \cdots & 0 \\
		-i\lambda \sigma_r & Y_r & 0 & \cdots & 0 \\
		-i\lambda\sigma_{r'} & Y_{r'} & 1 & \cdots & 0 \\
		\vdots & \vdots & \vdots & \ddots & \vdots \\
		-i\lambda\sigma_{r''} & Y_{r''} & 0 & \cdots & 1
	\EP,
\EE
where the sequence $(r',\ldots,r'')$ of length $m-1$ is the sequence $(1,2,\ldots,m)$ with the term $r$ excluded. Clearly, $\det \mathcal{B}_r = \det\mathcal{A}_r$. Hence
\BE
	\det\mathcal{A}_r = i\lambda\left( Y_r \sum_{\substack{p=1\\p\neq r}}^m \sigma_p - \sigma_r \sum_{\substack{p=1\\p\neq r}}^m Y_p \right),
\EE
and it follows that
\begin{align} \label{eqn:mSemiInf:Q0}
	g_0(\lambda^2,t) &= \frac{-\sum_{p=1}^m Y_p}{i\lambda\sum_{p=1}^m\sigma_p}, \\ \label{eqn:mSemiInf:Q1p}
	g_1^{r}(\lambda^2,t) &= \frac{Y_{r}\sum_{p=1}^m\sigma_p-\sigma_{r}\sum_{p=1}^m Y_p}{\sigma_r^2\sum_{p=1}^m\sigma_p}, & &1\leq r\leq m.
\end{align}
Substituting Equations~\eqref{eqn:mSemiInf:Q0} and~\eqref{eqn:mSemiInf:Q1p} into Equation~\eqref{eqn:ImplicitSol.Inf.t}, we obtain expressions for each $q^{r}$ in terms of the Fourier transforms of every initial datum $\hat{q}_0^{r}(\cdot)$ and  of every solution $\hat{q}^{r}(\cdot;t)$:
\begin{equation}\label{eqn:ExplicitSol.Inf.t}
\begin{split} 
	q^{r}(x,t) ={} &\frac{1}{2\pi} \int_{-\infty}^\infty e^{ix\lambda - \sigma_r^2\lambda^2t} \hat{q}_0^{r}(\lambda) \d \lambda \\
	&+ \int_{\partial D^+} \frac{e^{ix\lambda/\sigma_r-\lambda^2t}}{2\pi} \left( \frac{\sum_{p=1}^m\hat{q}_0^p\left(\frac{-\lambda}{\sigma_p}\right)}{\sum_{p=1}^m\sigma_p} - \frac{\hat{q}_0^r\left(\frac{-\lambda}{\sigma_r}\right)\sum_{p=1}^m\sigma_p - \sigma_r\sum_{p=1}^m\hat{q}_0^p\left(\frac{-\lambda}{\sigma_p}\right)}{\sigma_r\sum_{p=1}^m\sigma_p} \right) \d\lambda \\
	&- \int_{\partial D^+} \frac{e^{ix\lambda/\sigma_r}}{2\pi} \left( \frac{\sum_{p=1}^m\hat{q}^p\left(\frac{-\lambda}{\sigma_p};t\right)}{\sum_{p=1}^m\sigma_p} - \frac{\hat{q}^r\left(\frac{-\lambda}{\sigma_r};t\right)\sum_{p=1}^m\sigma_p - \sigma_r\sum_{p=1}^m\hat{q}^p\left(\frac{-\lambda}{\sigma_p};t\right)}{\sigma_r\sum_{p=1}^m\sigma_p} \right) \d\lambda.
\end{split}
\end{equation}
Note that we have separated the second integral of~\eqref{eqn:ImplicitSol.Inf.t} into two integrals, which is justified provided they converge. However, for $\lambda\to\infty$ from within the closure of $D^+$ ($\overline{D^+}$), $\Re(i\lambda)\to-\infty$ so by Jordan's Lemma the final integral of Equation~\eqref{eqn:ExplicitSol.Inf.t} not only converges but evaluates to $0$. Hence, we obtain an explicit integral representation of the solution,
\begin{equation}
\begin{split}
q^{r}(x,t) ={} &\frac{1}{2\pi} \int_{-\infty}^\infty e^{ix\lambda - \sigma_r^2\lambda^2t} \hat{q}_0^{r}(\lambda) \d \lambda \\
	&+ \int_{\partial D^+}\frac{ e^{ix\lambda/\sigma_r-\lambda^2t}}{2\pi} \left( \frac{2\sum_{p=1}^m\hat{q}_0^p\left(\frac{-\lambda}{\sigma_p}\right)}{\sum_{p=1}^m\sigma_p} - \frac{\hat{q}_0^r\left(\frac{-\lambda}{\sigma_r}\right)}{\sigma_r} \right) \d\lambda,
\end{split}
\end{equation}
in terms of only the initial data.

\begin{rmk}
In this example we were able to set $R=0$. This will not generally occur even for the semi-infinite rods when the network also contains finite rods (see, for example,~\cite[Proposition~3]{DeconinckPelloniSheils}). Indeed, $R=0$ is possible here because $\det\mathcal{A}$ is monomial which occurs for no configurations of rods except the one considered in this section.
\end{rmk}

\section{$m$ parallel finite rods} \label{sec:mParallel}
\begin{figure}
	\begin{center}
	\includegraphics{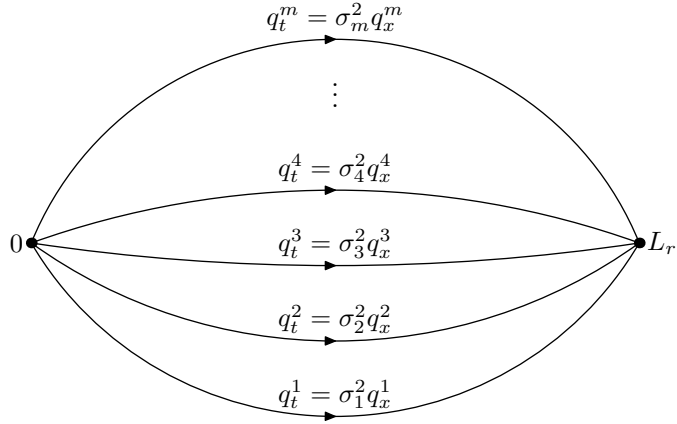}
	\caption{$m$ parallel finite rods.}
	\label{fig:mParallel}
	\end{center}
\end{figure}

In this section we consider the case of $m\in\N$ parallel finite rods of differing thermal diffusivity extending between a pair of vertices.  Note that we have chosen the rods to all be oriented one direction. This is purely for notational convenience, the parameterization of the rod from $0$ to $L_r$ or from $L_r$ to $0$ is not important.  The graph is as shown in Figure~\ref{fig:mParallel}. One of the two vertices is the terminus of every edge, and every edge emanates from the other vertex. As $|E|=m$, we find it notationally convenient to identify the set of edges $E$ with the set $\{1,2,\ldots,m\}$.

The continuity interface conditions allow us to define
\begin{subequations} \label{eqn:mParallel:IdentDirichlet}
\begin{align}
	g_0(\lambda^2,t) &:= g_0^{r}\left(\lambda^2,t\right), \\
	h_0(\lambda^2,t) &:= h_0^{r}\left(\lambda^2,t\right),
\end{align}
\end{subequations}
for $r=1,2,\ldots,m$.  Similarly, the continuity of flux interface conditions imply
\begin{subequations} \label{eqn:mParallel:IdentNeumann}
\begin{align}
	\sum_{r=1}^m \sigma_r^2 g_1^{r}(\lambda^2,t) &= 0, \\
	\sum_{r=1}^m \sigma_r^2 h_1^{r}(\lambda^2,t) &= 0.
\end{align}
\end{subequations}
Using the notation~\eqref{eqn:mParallel:IdentDirichlet} the global relation~\eqref{eqn:GR.Fin.t} reduces to
\begin{subequations}
\begin{equation}
\begin{split}\label{eqn:mParallel.GR.Fin}
	\sigma_r i\lambda g_0(\lambda^2,t) + \sigma_r^2 g_1^{r}(\lambda^2,t) - e^{\frac{-i\lambda L_r}{\sigma_r}} \left( \sigma_r i\lambda h_0(\lambda^2,t) + \sigma_r^2 h_1^{r}(\lambda^2,t) \right) \\
	= \hat{q}_0^{r}\left(\frac{\lambda}{\sigma_r}\right) - e^{\lambda^2t}\hat{q}^{r}\left(\frac{\lambda}{\sigma_r};t\right).
\end{split}
\end{equation}
Applying the transformation $\lambda\mapsto-\lambda$ which leaves the spectral functions $g_j^r(\lambda^2,t)$ and $h_j^r(\lambda^2,t)$ invariant we obtain
\begin{equation}
\begin{split}\label{eqn:mParallel.GR.Fin.n}
	- \sigma_r i\lambda g_0(\lambda^2,t) + \sigma_r^2 g_1^{r}(\lambda^2,t) - e^{\frac{i\lambda L_r}{\sigma_r}} \left( - \sigma_r i\lambda h_0(\lambda^2,t) + \sigma_r^2 h_1^{r}(\lambda^2,t) \right) \\
	= \hat{q}_0^{r}\left(\frac{-\lambda}{\sigma_r}\right) - e^{\lambda^2t}\hat{q}^{r}\left(\frac{-\lambda}{\sigma_r};t\right).
\end{split}
\end{equation}
\end{subequations}

This provides a system of $2m$ equations valid on $\C$ in $2m+2$ unknowns. System~\eqref{eqn:mParallel:IdentNeumann} provides two further equations in a subset of the same unknowns. We express this linear system as
$$
	\mathcal{A} X= Y,
$$
where
\begin{align} \notag
	X &= \left( g_0(\lambda^2,t), h_0(\lambda^2,t),\sigma_1^2 g_1^{1}(\lambda^2,t),\ldots, \sigma_m^2 g_1^{m}(\lambda^2,t), \sigma_1^2 h_1^{1}(\lambda^2,t),\ldots,\sigma_m^2 h_1^{m}(\lambda^2,t) \right)^\top, \\ \notag
	Y &= \left( 0,Y_1(\lambda),\ldots,Y_m(\lambda), Y_1(-\lambda),\ldots,Y_m(\lambda), \right)^\top, \\  \label{eqn:mParallel.Yp}
	Y_r(\lambda) &= \hat{q}_0^r\left(\frac{\lambda}{\sigma_r}\right) - e^{\lambda^2t}\hat{q}^r\left(\frac{\lambda}{\sigma_r};t\right),
\end{align}
and
\BE	\label{eqn:mParallel:A}
	\mathcal{A} =
	\BP
		0 & 0 & 1 & 1 & \cdots & 1 & 0 & 0 & \cdots & 0 \\
		0 & 0 & 0 & 0 & \cdots & 0 & 1 & 1 & \cdots & 1 \\
		i\lambda\sigma_1 & -i\lambda\sigma_1 e^{\frac{-i\lambda L_1}{\sigma_1}} &
			1 & 0 & \cdots & 0 & -e^{\frac{-i\lambda L_1}{\sigma_1}} & 0 & \cdots & 0 \\
		i\lambda\sigma_2 & -i\lambda\sigma_2 e^{\frac{-i\lambda L_2}{\sigma_2}} &
			0 & 1 & \cdots & 0 & 0 & -e^{\frac{-i\lambda L_2}{\sigma_2}} &  \cdots & 0 \\
		\vdots & \vdots & \vdots & \vdots & \ddots & \vdots & \vdots & \vdots & \ddots & \vdots \\
		i\lambda\sigma_m & -i\lambda\sigma_m e^{\frac{-i\lambda L_m}{\sigma_m}} &
			0 & 0 & \cdots & 1 & 0 & 0 & \cdots & -e^{\frac{-i\lambda L_m}{\sigma_m}} \\
		-i\lambda\sigma_1 & i\lambda\sigma_1 e^{\frac{i\lambda L_1}{\sigma_1}} &
			1 & 0 & \cdots & 0 & -e^{\frac{i\lambda L_1}{\sigma_1}} & 0 & \cdots & 0 \\
		-i\lambda\sigma_2 & i\lambda\sigma_2 e^{\frac{i\lambda L_2}{\sigma_2}} &
			0 & 1 & \cdots & 0 & 0 & -e^{\frac{i\lambda L_2}{\sigma_2}}  & \cdots & 0 \\
		\vdots & \vdots & \vdots & \vdots & \ddots & \vdots & \vdots & \vdots & \ddots & \vdots \\
		-i\lambda\sigma_m & i\lambda\sigma_m e^{\frac{i\lambda L_m}{\sigma_m}} &
			0 & 0 & \cdots & 1 & 0 & 0 & \cdots & -e^{\frac{i\lambda L_m}{\sigma_m}}
	\EP.
\EE

Solving the system via Cramer's rule, we obtain
\BE \label{eqn:mParallel.Q}
	\begin{aligned}
		g_0(\lambda^2,t) &= \frac{AC-BD}{\lambda(A^2-B^2)}, \qquad\qquad\qquad\qquad
		h_0(\lambda^2,t) = \frac{BC-AD}{\lambda(A^2-B^2)}, \\
		g_1^{r}(\lambda^2,t) &= \frac{-\det\BP B & A & C \\ A & B & D \\ \sigma_r\frac{1}{S_r(\lambda)} & \sigma_r\frac{C^{r}(\lambda)}{S_r(\lambda)} & \frac{{Y}_r(\lambda)e^{\frac{i\lambda L_r}{\sigma_r}}-{Y}_r(-\lambda)e^{\frac{-i\lambda L_r}{\sigma_r}}}{e^{\frac{i\lambda L_r}{\sigma_r}}-e^{\frac{-i\lambda L_r}{\sigma_r}}} \EP}{\sigma_r^2(A^2-B^2)}, \\
		h_1^{r}(\lambda^2,t) &= \frac{\det\BP B & A & D \\ A & B & C \\ \sigma_r\frac{1}{S_r(\lambda)} & \sigma_r\frac{C^{r}(\lambda)}{S_r(\lambda)} & \frac{{Y}_r(\lambda)-{Y}_r(-\lambda)}{e^{\frac{i\lambda L_r}{\sigma_r}}-e^{\frac{-i\lambda L_r}{\sigma_r}}} \EP}{\sigma_r^2(A^2-B^2)},
	\end{aligned}
\EE
where
\begin{align*}
	A &= \sum_{p=1}^m\sigma_p\frac{1}{S_p(\lambda)}, & C &= \sum_{p=1}^m\frac{Y_p(\lambda)-Y_p(-\lambda)}{e^{\frac{i\lambda L_p}{\sigma_p}}-e^{\frac{-i\lambda L_p}{\sigma_p}}}, \\
	B &= \sum_{p=1}^m\sigma_p\frac{C_p(\lambda)}{S_p(\lambda)}, & D &= \sum_{p=1}^m\frac{Y_p(\lambda)e^{\frac{i\lambda L_p}{\sigma_p}}-Y_p(-\lambda)e^{\frac{-i\lambda L_p}{\sigma_p}}}{e^{\frac{i\lambda L_p}{\sigma_p}}-e^{\frac{-i\lambda L_p}{\sigma_p}}},
\end{align*}
for
\begin{align*}
	S_p(\lambda) &= \sin(\lambda L_p/\sigma_p), & C_p(\lambda) &= \cos(\lambda L_p/\sigma_p).
\end{align*}

Equation~\eqref{eqn:mParallel.Q} for the $t$-transformed boundary and interface values depends upon the Fourier transform of the solution $\hat{q}^r(\cdot;t)$ through Expression~\eqref{eqn:mParallel.Yp} for $Y_p$. In order for Equation~\eqref{eqn:ImplicitSol.Fin.t} to represent an effective integral representation for each solution $q^r$ we must remove this dependence.

The product of a meromorphic function with a holomorphic function must include as zeros, every zero of the meromorphic function. As
\BE
	\left( \prod_{p=1}^{m} S_p(\lambda) \right)^2(A^2-B^2)
\EE
is an exponential polynomial (indeed an exponential sum), it is possible to obtain bounds on its zeros. Indeed, by~\cite[Theorem 3]{Lan1931a}, the zeros of the above product are confined within a finite-width horizontal strip. Hence, for sufficiently large $R>0$, there are no zeros of $A^2-B^2$ within $\overline{D_R^+}\cup\overline{D_R^-}$. Indeed, a simple calculation reveals
\BES
	R = \frac{(m-1)\log 2 + \log(8+m(m-1))}{\sqrt{2}\min_{1\leq r\leq m}\left\{\frac{L_r}{\sigma_r}\right\}}
\EES
is sufficient. A standard asymptotic analysis (see, for example,~\cite{DeconinckPelloniSheils}) shows that the terms involving $\hat{q}^p(\cdot;t)$ decay as $\lambda\to\infty$ from within the appropriate domains $D_R^\pm$, and Jordan's Lemma establishes that these terms make no contribution to the integral representation.

We have established the effective integral representation for the solution given by equation~\eqref{eqn:ImplicitSol.Fin.t} with spectral data given by equations~\eqref{eqn:mParallel.Q} where
$$
	Y_r(\lambda) = \hat{q}_0^r\left(\frac{\lambda}{\sigma_r}\right), \qquad 1\leq r\leq m.
$$

\section{$m$ finite rods connected at a single point} \label{sec:mFin}
\begin{figure}
	\begin{center}
	\includegraphics{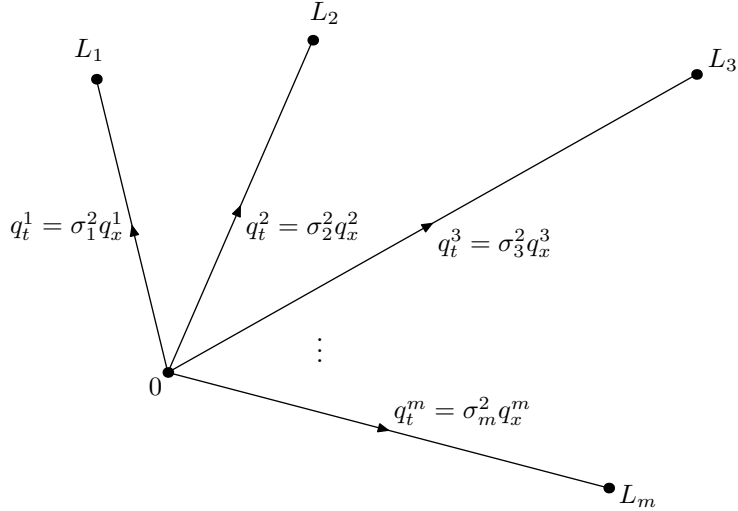}
	\caption{$m$ finite rods joined at a single vertex}
	\label{fig:mFin}
	\end{center}
\end{figure}

In this section we consider the case of $m\in\N$ finite rods of differing thermal diffusivity joined at a single vertex.  The graph is as shown in Figure~\ref{fig:mFin}. Of the $m+1$ vertices $m$ are the terminus of a single edge and all $m$ edges emanate from the other vertex. As $|E|=m$, we find it notationally convenient to identify the set of edges $E$ with the set $\{1,2,\ldots,m\}$. We enumerate the vertices so that each edge $r$ emanates from vertex $0$ and terminates at vertex $r$.

We use the condition of continuity across the interface to establish 
\BE \label{eqn:mfinite:D.Interface}
	g_0(\lambda^2,t):=g_0^{r}(\lambda^2,t)=\int_0^te^{\lambda^2s}\theta_0(s)\ud s
\EE
for $r=1,2,\ldots,m$, $\lambda\in\C$.

Recall the Robin boundary condition at $L_r$ for each of the terminal vertices with only one edge,~\eqref{RobinBC}. Taking the time transform of~\eqref{RobinBC} we have
\BE \label{RobinBC.t}
	\tilde{f}_r(\lambda^2;t) := \int_0^t e^{\lambda^2 s}f_r(s)\ud s = \beta_0^rh_0^{r}(\lambda^2,t)+\beta_1^rh_1^{r}(\lambda^2,t).
\EE
In order to reduce the size of the linear system that will eventually have to be solved, it is helpful to use the Robin condition to express either the 0\Th order or the 1\st order boundary value at $L_r$ in terms of the other. Suppose that precisely $1\leq m_N\leq m$ of the boundary conditions are in fact Neumann conditions (that is the corresponding $\beta^{r}_0=0$). Moreover suppose, without loss of generality, that these Neumann conditions are the boundary conditions at vertices $1,2,\ldots,m_N$, that $\beta_1^{r}=1$ for each $1\leq r\leq m_N$, and that $\beta_0^{r}=1$ for each $m_N+1\leq r\leq m$. We use the boundary conditions~\eqref{RobinBC.t} and Equation~\eqref{eqn:mfinite:D.Interface} to rewrite the global relation~\eqref{eqn:GR.Fin.t} as
\begin{subequations} \label{eqn:mfinite:GR}
\begin{equation}
\begin{split} \label{eqn:mfinite:GR.Neumann}
	\sigma_r i\lambda g_0(\lambda^2,t) + \sigma_r^2 g_1^{r}(\lambda^2,t) - e^{\frac{-i\lambda L_r}{\sigma_r}}i\lambda \sigma_r h_0^{r}(\lambda^2,t) = \\
	\hat{q}_0^{r}\left(\frac{\lambda}{\sigma_r}\right) - e^{\lambda^2t}\hat{q}^{r}\left(\frac{\lambda}{\sigma_r};t\right) + e^{\frac{-i\lambda L_r}{\sigma_r}}\sigma_r^2 \tilde{f}_r(\lambda^2;t),
\end{split}
\end{equation}
for $1\leq r\leq m_N$, $\lambda\in\C$ and
\begin{equation}
\begin{split} \label{eqn:mfinite:GR.Robin.Dirichlet}
	\sigma_r i\lambda g_0(\lambda^2,t) + \sigma_r^2 g_1^{r}(\lambda^2,t) - e^{\frac{-i\lambda L_r}{\sigma_r}} (-i\lambda\beta_1^r/\sigma_r+1) \sigma_r^2 h_1^{r}(\lambda^2,t) = \\
	\hat{q}_0^{r}\left(\frac{\lambda}{\sigma_r}\right) - e^{\lambda^2t}\hat{q}^{r}\left(\frac{\lambda}{\sigma_r};t\right) + e^{\frac{-i\lambda L_r}{\sigma_r}}\sigma_r i\lambda \tilde{f}_r(\lambda^2;t),
\end{split}
\end{equation}
\end{subequations}
or $m_N+1\leq r\leq m$ and $\lambda\in\C$.  Another set of $m$ global relations is obtained by evaluating equations~\eqref{eqn:mfinite:GR} under the map $\lambda\mapsto-\lambda$, and observing that the spectral functions $g_j$, $h_j$ and the $t$-transformed data $\tilde{f}_r$ are invariant under this map.

Together with a $t$-transform of the interface condition~\eqref{cont_qx} the $2m$ global relations form a system of linear equations in the $2m+1$ spectral functions
\begin{align*}
	&g_0(\lambda^2, t),&&&t>0, \\
	&g_1^{r}(\lambda^2, t),   & &1\leq r \leq m, &t>0, \\
	&h_0^{r}(\lambda^2, t), & &1\leq r\leq m_N,&t>0, \\
	&h_1^{r}(\lambda^2, t), &  &m_N+1\leq r\leq m,&t>0.
\end{align*}
The system can be expressed as
\begin{subequations} \label{eqn:mfinite:Linear.System}
\BE
	\mathcal{A}X = Y,
\EE
where
\begin{align}
	X &= \big( g_0, \sigma_1^2 g_1^{1},\ldots,\sigma_m^2 g_1^{m}, \sigma_1 h_0^{1},\ldots,\sigma_{m_N} h_0^{m_N}, \sigma_{m_N+1}^2 h_1^{m_N+1},\ldots, \sigma_{m}^2 h_1^{m}\big)^\top , \\
	Y &= \left( 0, Y_1(k),\ldots,Y_m(k), Y_1(-k),\ldots,Y_m(-k) \right)^\top,\\
	Y_r(\lambda) &=
	\begin{cases}
		\hat{q}_0^r\left(\frac{\lambda}{\sigma_r}\right) - e^{\lambda^2t}\hat{q}^r\left(\frac{\lambda}{\sigma_r};t\right) + e^{\frac{-i\lambda L_r}{\sigma_r}}\sigma_r^2\tilde{f}_r(\lambda;t) & 1\leq r\leq m_N, \\
		\hat{q}_0^r\left(\frac{\lambda}{\sigma_r}\right) - e^{\lambda^2t}\hat{q}^r\left(\frac{\lambda}{\sigma_r};t\right) + i\lambda e^{\frac{-i\lambda L_r}{\sigma_r}}\sigma_r\tilde{f}_r(\lambda;t) &m_N+1\leq r\leq m,
	\end{cases} \\
	\mathcal{A} &=
	\BP
		0                       & 1_{1 \times m} & 0_{1\times m}           \\
		\mathcal{A}_0(\lambda)  & I_{m \times m} & \mathcal{A}_1(-\lambda) \\
		\mathcal{A}_0(-\lambda) & I_{m \times m} & \mathcal{A}_1(\lambda)  \\
	\EP, \\
	\mathcal{A}_0(\lambda) &= \big( i\lambda \sigma_1, \ldots, i\lambda \sigma_m \big)^\top,
\end{align}
and
\begin{multline}
	\mathcal{A}_1(\lambda) = \\
	\BP
		i\lambda e^{\frac{i\lambda L_1}{\sigma_1}} & \cdots & 0 & 0 & \cdots & 0 \\
		\vdots & \ddots & \vdots & \vdots & \ddots & \vdots \\
		0 & \cdots & i\lambda e^{\frac{i\lambda L_{m_N}}{\sigma_{m_N}}} & 0 & \cdots & 0 \\
		0 & \cdots & 0 & \left(\frac{-i\lambda\beta_1^{m_N+1}}{\sigma_{m_N+1}}-1\right)e^{\frac{i\lambda L_{m_N+1}}{\sigma_{m_N+1}}} & \cdots & 0 \\
		\vdots & \ddots & \vdots & \vdots & \ddots & \vdots \\
		0 & \cdots & 0 & 0 & \cdots & \left(\frac{-i\lambda\beta_1^{m}}{\sigma_m}-1\right)e^{\frac{i\lambda L_m}{\sigma_m}}
	\EP.
\end{multline}
\end{subequations}
The boundary conditions~\eqref{RobinBC.t} give the remaining spectral functions in terms of those appearing in the vector $X$.

Solving this system and substituting the solution into equation~\eqref{eqn:ImplicitSol.Fin.t} provides an integral representation of the solution in terms of the initial and boundary data, and the Fourier transform of the solution $\hat{q}^r(\cdot;t)$. It remains to show that it is possible to remove the dependence on the Fourier transform of the solution.

When the system~\eqref{eqn:mfinite:Linear.System} is solved using Cramer's rule, each entry $X_j$ of $X$ is expressed as the ratio of two determinants, each of which is an analytic function of $\lambda$. The denominator is an exponential polynomial in which every exponent is purely imaginary. Hence the zeros of the denominator lie asymptotically within a horizontal logarithmic strip~\cite[Theorem~6]{Lan1931a}. Moreover, $X_j$ is holomorphic except at the zeros of the denominator, so choosing $R$ sufficiently large ensures that $X_j$ is holomorphic on $\overline{D^+_R}\cup\overline{D^-_R}$.

As the growing (respectively, bounded) entries of $\mathcal{A}$ lie in the same rows as the growing (respectively, bounded) entries of $Y$ involving $\hat{q}^r(\cdot;t)$, but each such entry in $Y$ grows (respectively, decays) like $\lambda^{-1}$ multiplied by the corresponding entry in $\mathcal{A}$, each integrand involving $\hat{q}^r(\cdot;t)$ decays as $\lambda\to\infty$ from within the relevant sector $\overline{D^\pm_R}$. Hence, by Jordan's Lemma, the terms involving $\hat{q}^r(\cdot;t)$ make no contribution to the solution representation.

\begin{rmk}
Using a similar Jordan's Lemma argument, it can be justified that $g_j^r(\lambda^2,T)$ and $h_j^r(\lambda^2,T)$ may be used in place of $g_j^r(\lambda^2,t)$ and $h_j^r(\lambda^2,t)$
(see, for example,~\cite[page~8]{Fok2008a}).
Evaluating these $t$-transforms at final time $T$ only is helpful for numerical calculations.
\end{rmk}

\subsection{Example}

We study the problem of three finite rods emanating from a single vertex. The rods have lengths
\BE
	L_1 = 1, \qquad L_2 = 1, \qquad L_3 = 2.
\EE
Dirichlet conditions are prescribed at the terminal vertices
\BE
	q^1(1,t) = \sin(t), \qquad q^2(1,t) = 0, \qquad q^3(2,t) = 0,
\EE
so that the system is driven by an oscillating temperature profile at the terminal vertex of rod 1. The physical properties of the rods are such that their scaled thermal diffusivities are
\BE
	\sigma_1^2 = 4, \qquad \sigma_2^2 = 9, \qquad \sigma_3^2 = 1.
\EE

For this system, equation~\eqref{eqn:mfinite:Linear.System} evaluates to
\BE \label{eqn:mfinite.Eg:Linear.System}
	\BP
		0 & 1 & 1 & 1 & 0 & 0 & 0 \\
		 2i\lambda & 1 & 0 & 0 & - e^{ i\lambda/2} & 0 & 0 \\
		 3i\lambda & 0 & 1 & 0 & 0 & - e^{ i\lambda/3} & 0 \\
		  i\lambda & 0 & 0 & 1 & 0 & 0 & - e^{ 2i\lambda } \\
		-2i\lambda & 1 & 0 & 0 & - e^{-i\lambda/2} & 0 & 0 \\
		-3i\lambda & 0 & 1 & 0 & 0 & - e^{-i\lambda/3} & 0 \\
		 -i\lambda & 0 & 0 & 1 & 0 & 0 & - e^{-2i\lambda }
	\EP
	\BP
		g_0 \\ 4g_1^1 \\ 9g_1^2 \\ g_1^3 \\ 4h_1^1 \\ 9h_1^2 \\ h_1^3
	\EP
	= Y.
\EE
The determinant of this system is
\BE
	-4i\lambda \left( 2\cos\left(\frac{7\lambda}{6}\right) + \cos\left(\frac{11\lambda}{6}\right) - 3\cos\left(\frac{17\lambda}{6}\right) \right),
\EE
which has no nonreal zeros so, for any choice of $R>0$, the terms involving $\hat{q}^r(\cdot,t)$ make no contribution to the solution representation. Therefore, for the purposes of numerical work, it is sufficient to solve the system~\eqref{eqn:mfinite.Eg:Linear.System} with
\BE
	Y = \left( 0, 2i\lambda \int_0^T e^{\lambda^2 s} \sin(s) \ud\lambda, 0, 0, -2i\lambda \int_0^T e^{\lambda^2 s} \sin(s) \ud\lambda, 0, 0 \right)^\top
\EE
and substitute the resulting formulae for $g$, $h$ into equation~\eqref{eqn:ImplicitSol.Derivation.1.Fin}.

\begin{figure}
	\begin{center}
	\begin{minipage}[c]{10.2cm}\includegraphics[width=10cm]{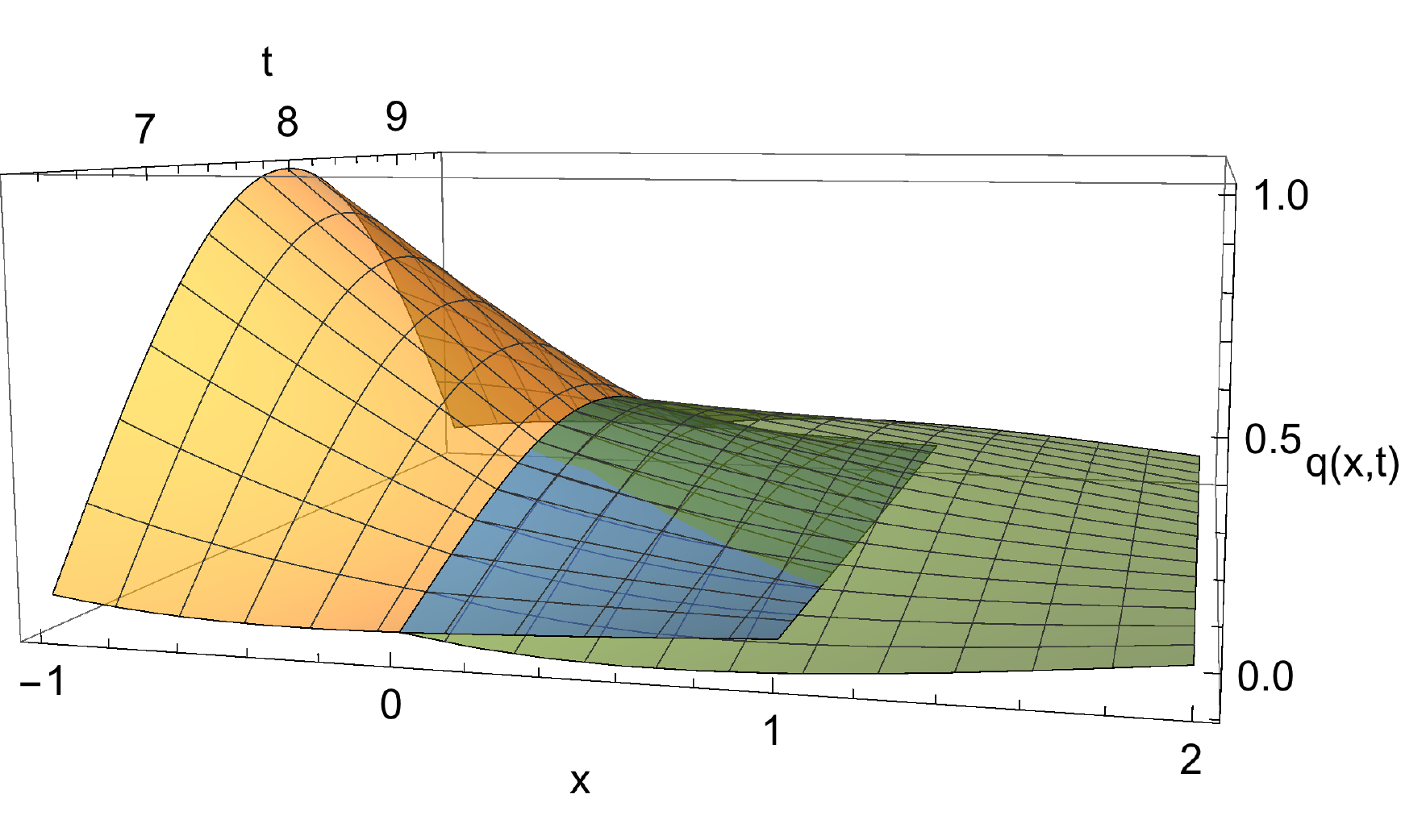}\end{minipage}\begin{minipage}[c]{2.1cm}\includegraphics[width=2cm]{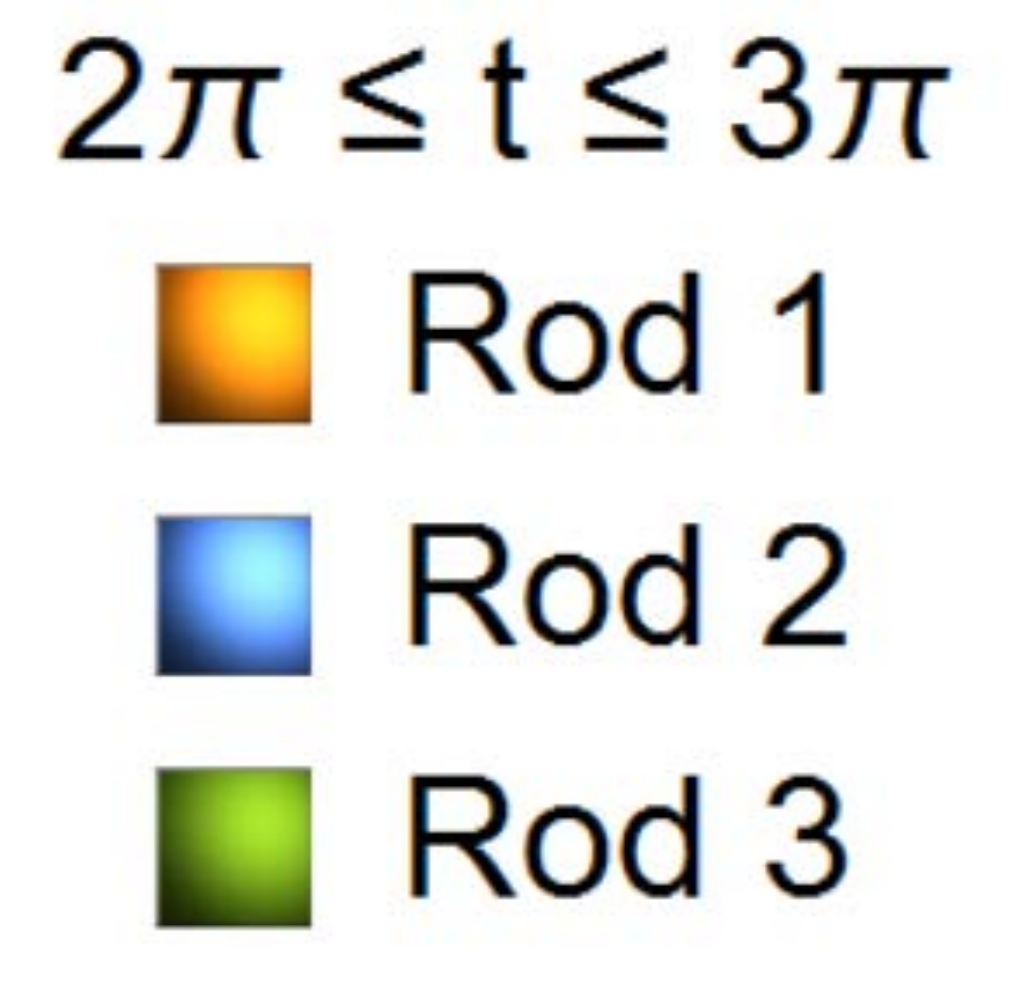}\end{minipage}
	\label{fig:3finrod1}
	\end{center}
	\begin{center}
	\begin{minipage}[c]{10.2cm}\includegraphics[width=10cm]{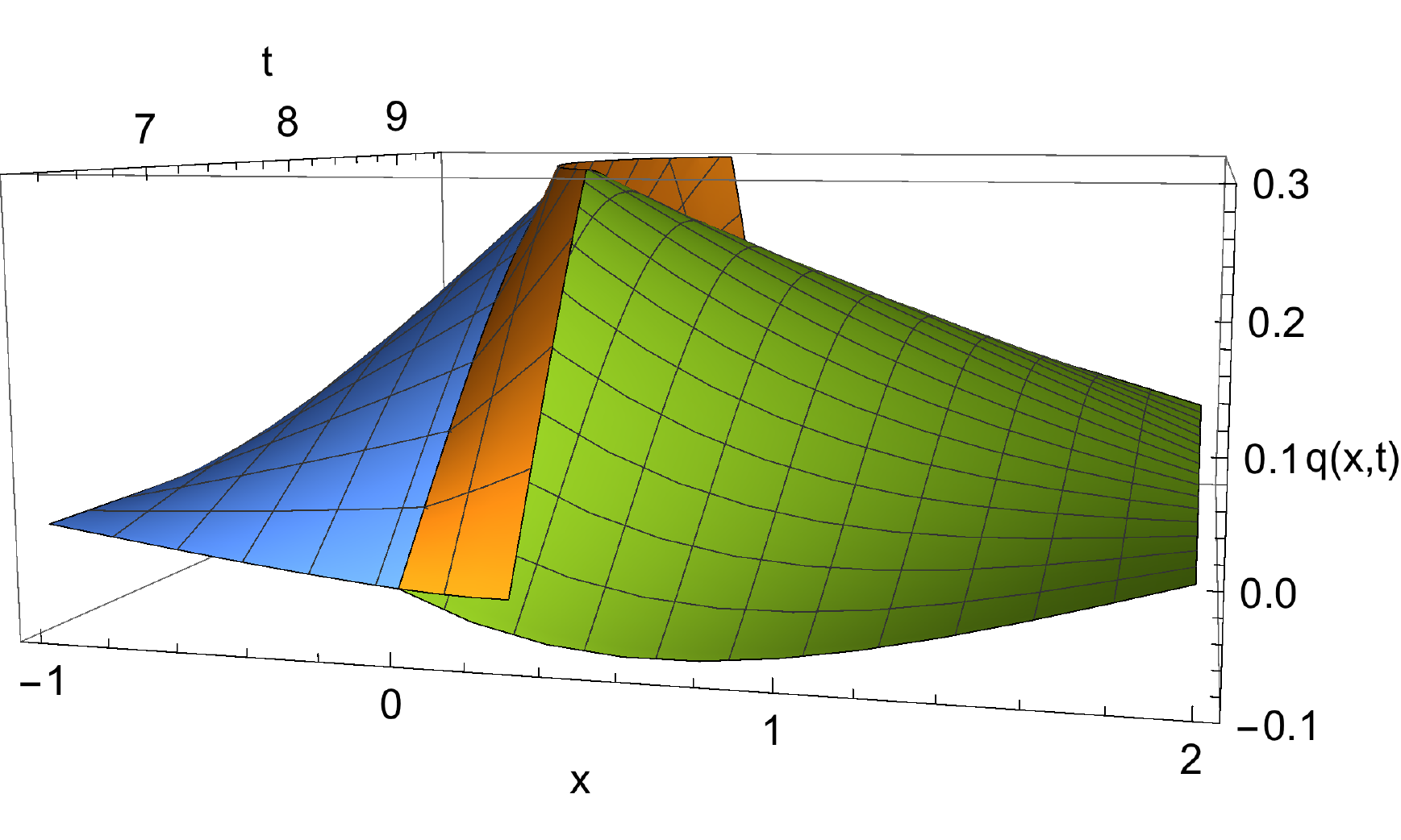}\end{minipage}\begin{minipage}[c]{2.9cm}\includegraphics[width=2.8cm]{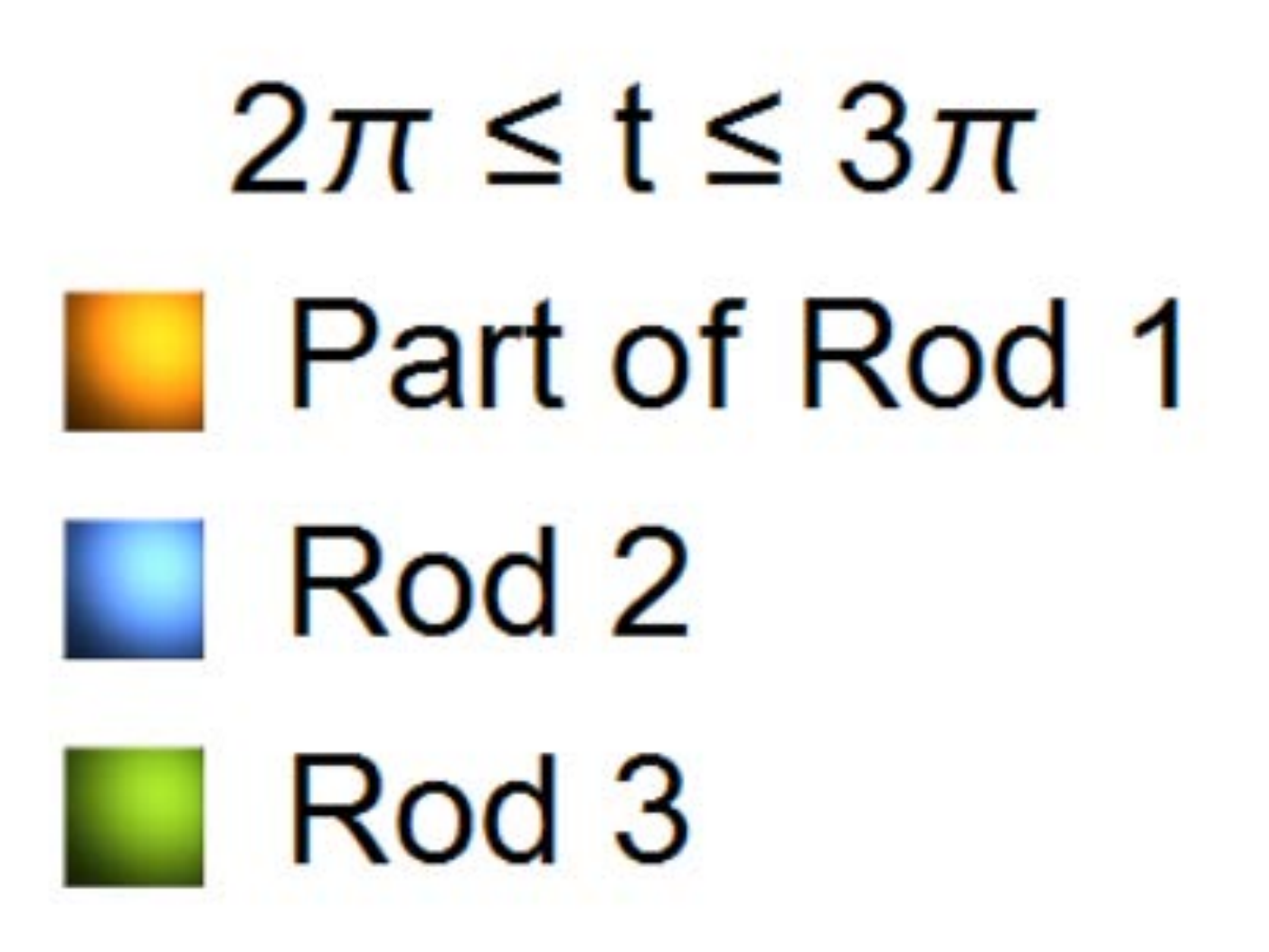}\end{minipage}
	\caption{Three finite rods with $2\pi \leq t \leq 3\pi$. The heat distributions $q^1(x,t)$ (orange), $q^2(x,t)$ (blue), and $q^3(x,t)$ (green) are shown. The code used to produce this plot was written by R.\ J.\ Buckingham and D.\ A.\ Smith.}
	\label{fig:3finrod2}
	\end{center}
\end{figure}

\begin{figure}
	\begin{center}
	\includegraphics{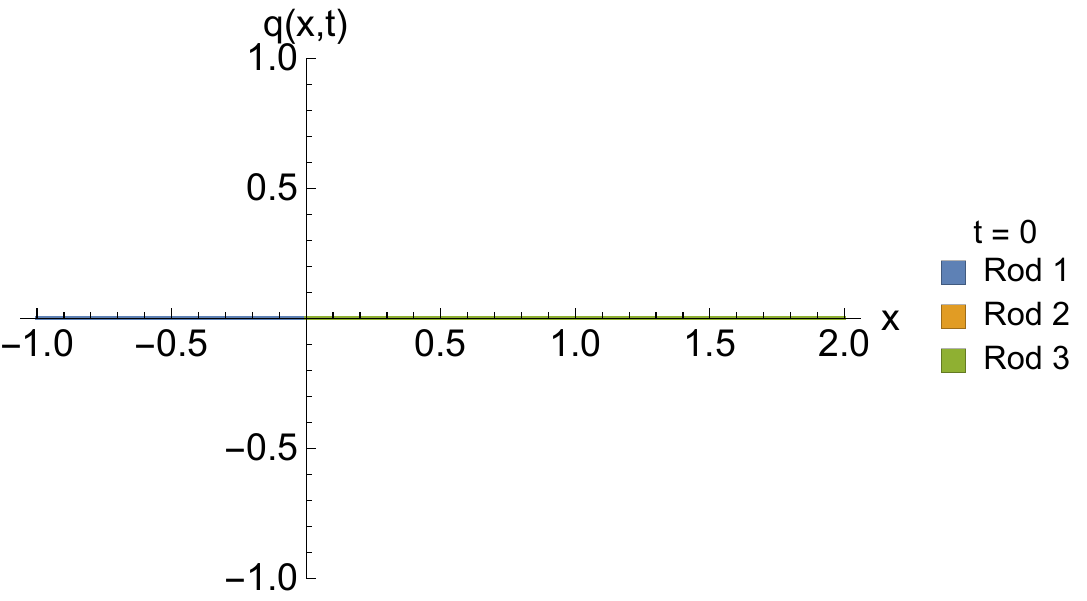}
	\caption{{\it Please see the ancillary file ``threerodsplotfrom0-oneshot.swf''.} Three finite rods with $0 \leq t \leq 3\pi$. The heat distributions $q^1(x,t)$ (orange), $q^2(x,t)$ (blue), and $q^3(x,t)$ (green) are shown. The code used to produce this plot was written by R.\ J.\ Buckingham and D.\ A.\ Smith.}
	\label{fig:3finrod3}
	\end{center}
\end{figure}

\begin{figure}
\centering
	\includegraphics[width=.45\textwidth]{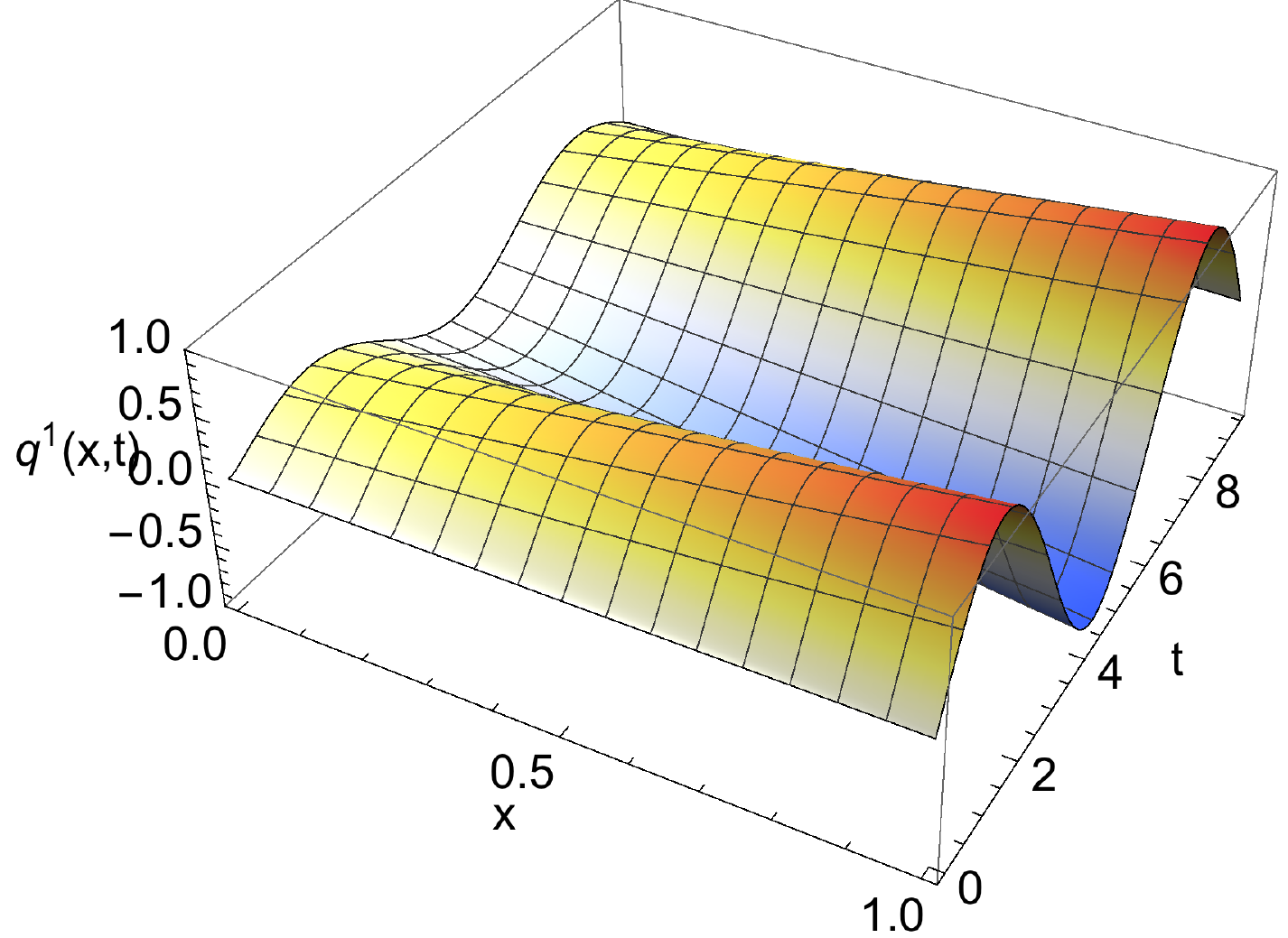}
		\includegraphics[width=.45\textwidth]{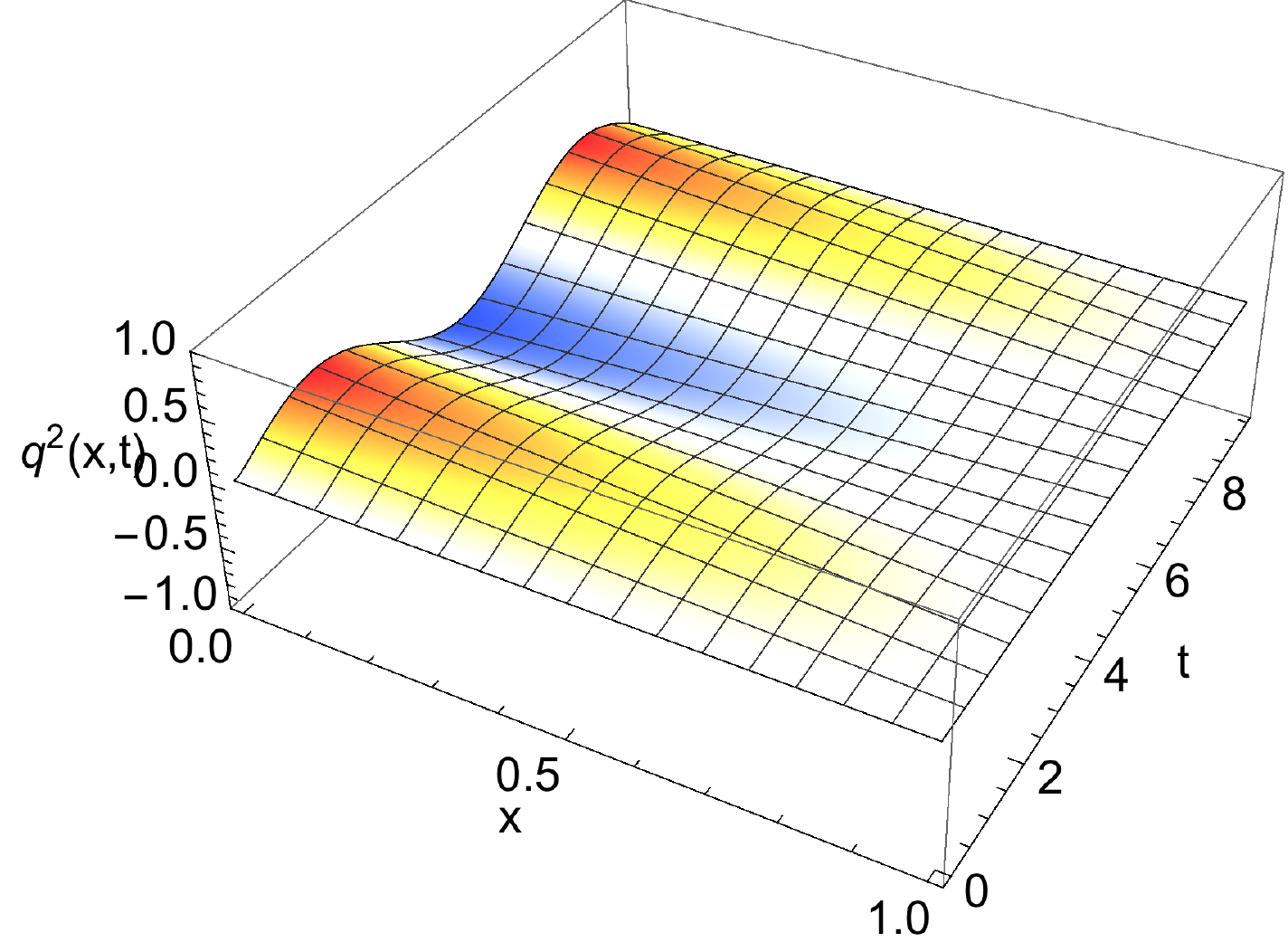}
			\includegraphics[width=.45\textwidth]{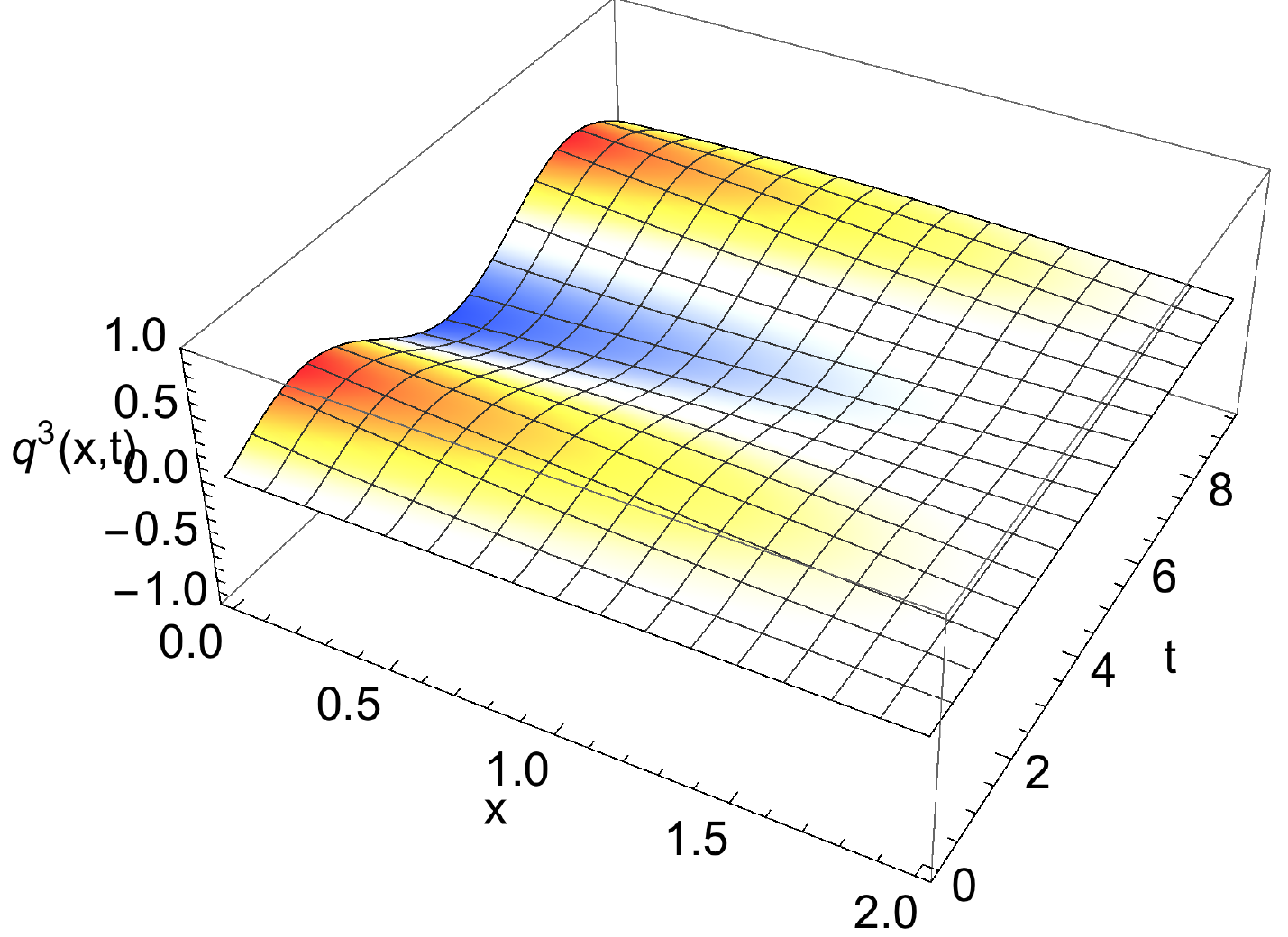}
\caption{Three finite rods $q^1(x,t)$, $q^2(x,t)$ and $q^3(x,t)$ with $0\leq t \leq 3\pi$. The code used to produce this plot was written by R.\ J.\ Buckingham and D.\ A.\ Smith.}
	\label{fig:3finrod3}
\end{figure}

Using a Mathematica code written by R.\ J.\ Buckingham and D.\ A.\ Smith, we obtain the plots of $q^1$ (orange), $q^2$ (blue), and $q^3$ (green), for $x\in[0,L_r]$ and $t\in[0,2\pi]$ as shown in Figures~\ref{fig:3finrod1} and~\ref{fig:3finrod2}. An animation is shown in Figure~\ref{fig:3finrod3} showing all $q^r$, for $0 \leq t \leq 3\pi$.

\section{Concluding remarks}

\begin{rmk}
In this work we have provided a general framework implementing the Fokas method for solving interface problems for the heat equation on a network. We have demonstrated the full argument in a number of cases, but it should be noted that the method is generally applicable to any such problem. A general argument for the effectiveness of the generalized spectral Dirichlet-to-Neumann map, including the Jordan's Lemma argument that terms involving $\hat{q}^r(\lambda;t)$ do not contribute to the resulting solution representation, was omitted in this paper but will appear elsewhere.
\end{rmk}
\begin{rmk}
The examples presented in this paper were selected because their generalized spectral Dirichlet-to-Neumann maps can be formulated as linear systems which are relatively easy to solve by hand. It should be noted that this is not always the case. In general, for a graph with $m_f$ finite edges and $m_i$ infinite edges, it will always be possible to express the generalized spectral Dirichlet-to-Neumann map as a linear system of dimension $4m_f+2m_i$ but it may not be practical to solve this system by hand. For example, the case of $m$ finite rods joined end-to-end in a line (see~\cite{DeconinckPelloniSheils} for the cases $m=2$ and $m=3$, see~\cite{Asvestas} for arbitrary $m$ with homogeneous Dirichlet conditions) has a generalized spectral Dirichlet-to-Neumann map expressible as a pentadiagonal linear system. It is the authors' opinion that a computer code would be the most appropriate tool to obtain the solution representation for arbitrarily complex networks.
\end{rmk}
\begin{rmk}
It has been established that the Fokas method may be used to solve all well-posed initial-boundary value problems for a linear constant-coefficient evolution equation on a finite interval or half-line~\cite{Fok2000a,FP2001a,Smi2012a}. A number of other interface problems have already been studied using the Fokas method for the heat equation~\cite{Asvestas,DeconinckPelloniSheils,SheilsDeconinck_PeriodicHeat,Mantzavinos}, the linear Schr\"{o}dinger equation~\cite{SheilsDeconinck_LS, SheilsDeconinck_LSp}, and a third-order dispersive equation~\cite{DeconinckSheilsSmith}.  Given the results above, it is reasonable to expect that the method could be used to solve problems for higher order PDE on a general network.

In the present work, we made use of a particular physical problem to determine appropriate interface conditions but such a physical model may not always be available. However, as the Fokas method can be used to determine well-posedness criteria for half-line and finite interval problems~\cite{Pel2004a,Smi2011a,Smi2012a,Smi2013a}, it should be possible to determine interface and boundary conditions that specify well-posed problems for higher order equations on networks.
\end{rmk}

\section*{Acknowledgements}
This work was generously supported by the National Science Foundation under grant number NSF-DGE-0718124 (N.E.S.).  The authors wish to thank Bernard Deconinck for helpful discussions.  Any opinions, findings, and conclusions or recommendations expressed in this material are those of the authors and do not necessarily reflect the views of the funding sources.

\bibliographystyle{amsplain}
\bibliography{dbrefs,FullBib}
\end{document}